\documentclass[oneside]{amsart}
\usepackage{amsmath, amssymb}
\usepackage{amsthm}
\newtheorem{theorem}{Theorem}[section]
\newtheorem{proposition}[theorem]{Proposition}
\newtheorem{lemma}[theorem]{Lemma}
\newtheorem{corollary}[theorem]{Corollary}

\newtheorem{example}{Example}

\def\F{\mathcal{F} }

\def\S{\mathbb{S} } 

\def\T{\mathbb{T} }

\def\R{\mathbb{R} } 
\def\Z{\mathbb{Z} } 
\def\nbd{neighborhood } 
 
\def\R{\mathbb{R} }

\title[A topological characterization for non-wandering surface flows]{A topological characterization for non-wandering surface flows}

\author{Tomoo YOKOYAMA}

\address{ Department of Mathematics, Faculty of Education, Kyoto University of Education, 
              1 Fujinomori, Fukakusa, Fushimi-ku, Kyoto, 612-8522, Japan} 
\email{tomoo@kyokyo-u.ac.jp}

\date{\today}

\thanks{The author is partially supported 
by the JST CREST Program at Department of Mathematics,  
Kyoto University of Education.}

\begin{document}

\maketitle

\begin{abstract}
Let $v$ be a continuous flow with arbitrary singularities 
on a compact surface. 
Then 
we show that 
if $v$ is non-wandering 
then 
$v$ is topologically equivalent to a $C^{\infty}$ flow 
such that  
there are no exceptional orbits 
and 
$\mathrm{P} \sqcup \mathop{\mathrm{Sing}}(v) = 
\{ x \in M \mid \omega(x) \cup \alpha(x) \subseteq \mathop{\mathrm{Sing}}(v) \}$,  
where $\mathrm{P}$ is the union of non-closed proper orbits 
and $\sqcup$ is the disjoint union symbol.  
Moreover, 
$v$ is non-wandering 
if and only if 
$\overline{\mathrm{LD}\sqcup \mathop{\mathrm{Per}}(v)}  
\supseteq M - \mathop{\mathrm{Sing}}(v)$,  
where 
$\mathrm{LD}$ is the union of locally dense orbits 
and $\overline{A}$ is the closure of a subset $A \subseteq M$.  
On the other hand, 
$v$ is topologically transitive 
if and only if 
$v$ is non-wandering 
such that 
$ \mathop{\mathrm{int}}(\mathop{\mathrm{Per}}(v) 
\sqcup \mathop{\mathrm{Sing}}(v)) = \emptyset$ 
and 
$M - (\mathrm{P}  \sqcup \mathop{\mathrm{Sing}}(v))$ is connected, 
where $\mathrm{int} {A}$ is the interior of a subset $A \subseteq M$.  
In addition, 
we construct 
a smooth flow on $\mathbb{T}^2$ 
with 
$\overline{\mathrm{P}} = \overline{\mathrm{LD}} =\mathbb{T}^2$. 
\end{abstract}


\section{Introduction and preliminaries}

In  \cite{P}, H. Poincar\'e has constructed a flow on a torus 
with an exceptional minimal set. 
A. Denjoy has constructed such a $C^1$-flow \cite{D}. 
On the other hand,  
the author has shown that 
there are no exceptional minimal sets of $C^2$-flows on tori. 
This result is generalized to the compact surface cases by Schwartz, A. J. \cite{S}.
%
%
%
On the other hand, 
in \cite{NZ}, 
they have shown the characterization of 
the non-wandering flows on compact surfaces with finitely many singularities. 
In \cite{Mar}, 
one has given a description near orbits of the non-wandering flow with 
the set of singularities which is totally disconnected. 
Moreover, 
in  \cite{M}, 
it has shown that 
if a nontrivial recurrent point $x$ of a surface flow 
belongs to the limit
set of another nontrivial recurrent point $y$, 
then 
$y$ 
belongs to the limit set of $x$. 
We sharpen these results to analysis surface flows. 
In particular, 
it has shown that 
the orbit class of each 
nontrivial weakly recurrent point is 
the orbit closure in 
the set of regular weakly recurrent points. 
In this paper, 
we show the non-existence of exceptional minimal sets of 
continuous non-wandering flows 
with arbitrary singularities 
on compact surfaces 
and obtain  
a topological characterization of 
non-wandering flows 
and 
a topological characterization of 
non-periodic proper orbits 
of non-wandering flows.  
Moreover, we give 
a smoothability of non-wandering flows 
and 
a characterization of topological transitivity 
for continuous flows. 
In addition, 
we construct a smooth flow on $\mathbb{T}^2$ 
and a foliation on an open manifold contained in $\mathbb{R}^2$ 
such that 
$\mathrm{LD}$ and $\mathrm{P}$ are dense, 
where 
$\mathrm{LD}$ is the union of locally dense orbits (resp. leaves) 
and 
$\mathrm{P}$ is the union of non-closed proper orbits (resp. leaves). 
On the other hand, 
if $\mathop{\mathrm{Sing}}(v)$ 
consists of finitely many contractible connected components, 
then $\mathrm{LD}$ is open 
and $\mathrm{P}$ 
consists of finitely many orbits.

By flows, 
we mean continuous $\mathbb{R}$-actions on surfaces.
Let 
$v: \mathbb{R} \times M \to M$ be a flow on a compact surface $M$. 
Put $v_t( \cdot ) := v(t, \cdot )$ 
and $O_v( \cdot ) := v(\mathbb{R}, \cdot )$.  
A subset of $M$ is said to be saturated 
if 
it is a union of orbits. 
Recall that 
a point $x$ of $M$ is 
singular if 
$x = v_t(x)$ for any $t \in \R$, 
is regular if 
$x$ is not singular, 
and 
is periodic if 
there is positive number $T > 0$
such that 
$x = v_T(x)$ and  
$x \neq v_t(x)$ for any $t \in (0, T)$. 
Denote by 
$\mathop{\mathrm{Sing}}(v)$ 
(resp. $\mathop{\mathrm{Per}}(v)$) 
the set of singular (resp. periodic) points. 
A point $x$ is wandering if 
there are a neighbourhood $U$ of $x$ and 
a positive number $N$ 
such that 
$\bigcup_{t > N} v_t(U) \cap U = \emptyset$, 
and is  
non-wandering if  
$x$ is not wandering 
(i.e. for each neighbourhood $U$ of $x$ and 
each positive number $N$, 
there is $t \in \mathbb{R}$ with $|t| > N$ such that 
$v_t(U) \cap U \neq \emptyset$). 
An orbit is non-wandering if 
it consists of non-wandering points 
and 
the flow $v$ is non-wandering if 
every point is non-wandering. 
For a point $x \in M$, 
define 
the omega limit set $\omega(x)$
and 
the alpha limit set $\alpha(x)$ 
of $x$
as follows: 
$\omega(x) 
:= \bigcap_{n\in \mathbb{R}}\overline{\{v_t(x) \mid t > n\}} 
$, 
$\alpha(x) 
:= \bigcap_{n\in \mathbb{R}}\overline{\{v_t(x) \mid t < n\}} 
$. 
A point $x$ of $M$ is 
positive recurrent (resp. negative recurrent) if 
$x \in \omega(x)$ 
(resp.  $x \in \alpha(x)$),  
and that 
$x$ is recurrent (resp. weakly recurrent) 
if $x$ is positive and (resp. or) negative recurrent. 
A (weakly) recurrent orbit is an orbit of such a point. 
An orbit is proper if 
it is embedded, 
locally dense if 
the closure of it has nonempty interior, 
and 
exceptional if 
it is neither proper nor locally dense.  
A point is proper (resp. locally dense, exceptional) if 
so is its orbit. 
Denote by 
$\mathrm{LD}$ 
(resp. $\mathrm{E}$, $\mathrm{P}$)
the union of locally dense orbits 
(resp. exceptional orbits, 
non-closed proper orbits).  
Note 
$\mathrm{P}$ is the complement of 
the set of weakly recurrent points. 
%
By the definitions,  
we have a decomposition 
$\mathop{\mathrm{Sing}}(v) \sqcup 
\mathop{\mathrm{Per}}(v) \sqcup \mathrm{P}
\sqcup \mathrm{LD} 
\sqcup \mathrm{E} = M$. 
A (weakly) recurrent orbit is nontrivial if 
it is not closed.  
Note that 
the union of nontrivial weakly recurrent orbits 
corresponds with $\mathrm{LD} \sqcup \mathrm{E}$.   
A quasi-minimal set of $v$ is the closure of a nontrivial weakly recurrent orbit. 
It's known that 
the total number of quasi-minimal sets for $v$ 
cannot exceed $g$ if $M$ is an orientable surface of genus $g$ \cite{M}, 
and $\frac{p-1}{2}$ if $M$ is a non-orientable surface of genus $p$ \cite{Ma}. 
Therefore 
the closure 
$\overline{\mathrm{LD \sqcup E}}$ consists of finitely many quasi-minimal sets. 
By a limit cycle, we mean 
a periodic orbit of $v$ which is the $\alpha$-limit set 
or the $\omega$-limit set of some point 
not on the periodic orbit. 

\section{A Topological characterization of 
non-wandering surface flows with arbitrary singularities}

Let $v$ be a continuous flow on a compact surface $M$. 
We call 
that a collar $A$ of an periodic orbit $O \subset A$  is 
an annulus $A$ one of whose connected component of 
$\partial A$ is $O$, 
where $\partial A := \overline{A} - \mathrm{int} A$ 
is the topological boundary of $A$. 
First, we state a following easy observation. 

\begin{lemma}\label{lem00}
$O \subset \mathop{\mathrm{Per}}(v) \sqcup \mathrm{P}$ 
for an orbit $O$ with $\overline{O} \cap \mathop{\mathrm{Per}}(v) \neq \emptyset$.  
Moreover 
each limit cycle is contained in 
$\partial (\mathrm{int} \mathrm{P})$. 
%
\end{lemma}

\begin{proof}
Let $O$ be an non-periodic orbit with $\overline{O} \cap \mathop{\mathrm{Per}}(v) \neq \emptyset$. 
Then $\overline{O} - O$ contains a limit cycle $\gamma$. 
The  flow box theorem (cf. Theorem 1.1, p.45\cite{ABZ}) implies that 
a limit cycle $\gamma$ is covered by finitely many flow boxes $\{ U_i' \}$. 
Since $\overline{O} \cap \gamma \neq \emptyset$,  
the one-sided holonomy of $\gamma$ is contracting or expanding. 
Fix a point $z$ of $O$ in a flow box $U_i'$. 
Then the point $z$ has an open neighbourhood $U \subset U_i'$ 
which is a flow box such that 
$O \cap U$ is one arc. 
Then $\cup_{t \in \mathbb{R}} v_t(U)$ 
is an open neighbourhood of $O$ in 
which $O$ is closed.    
Hence $O \subset \mathrm{P}$.  
Let $U' \subset \bigcup_i U_i'$ be 
a sufficiently small collar of $\gamma$ 
where the holonomy along $\gamma$ is contracting or expanding.  
Then the orbit closure of each point $y$ in $U' - \gamma$ contains $\gamma$ 
but the orbit of $y$ is not closed. 
Therefore $y \in \mathrm{P}$
 and so $V := \cup_{t \in \mathbb{R}} v_t(U') - \gamma \subseteq \mathrm{P}$ 
is a saturated open subset with $\gamma \subseteq \partial V 
\subseteq \partial (\mathrm{int} \mathrm{P})$. 
This implies the second assertion. 
\end{proof}

Recall that 
the orbit class $\hat{O}$ of an orbit $O$ is 
the union of points each of whose orbit closure 
corresponds with $\overline{O}$ 
(i.e.  $\hat{O} := \{ y \in M \mid \overline{O} = \overline{O_v(y)} \}$). 
Now we show that 
the orbit class of each 
nontrivial weakly recurrent point is 
the orbit closure in 
the set of regular weakly recurrent points, 
which 
refine a Ma\v \i er-type result \cite{M}.  

\begin{proposition}\label{lem001}
For an 
orbit $O \subset \mathrm{LD} \sqcup \mathrm{E}$,  
we have 
$\hat{O} = \overline{O} 
\setminus (\mathop{\mathrm{Sing}}(v) \sqcup \mathrm{P}) 
\subset \mathrm{LD} \sqcup \mathrm{E}$. 
\end{proposition}

\begin{proof}
Let $Q$ be a quasi-minimal set which is the closure of $O$. 
By Lemma \ref{lem00}, 
we have $Q \cap \mathop{\mathrm{Per}}(v)  = \emptyset$. 
Note that 
the inverse image of 
proper (resp. locally dense, exceptional) orbits 
by any finite covering 
are also proper (resp. locally dense, exceptional). 
By taking a double covering of $M$ and the doubling of $M$ 
if necessary, 
we may assume that $v$ is transversally orientable 
and $M$ is closed and orientable. 
For any point $y \in \hat{O} - O$, 
we have 
$\overline{O_v(y)} = \overline{O}$ 
and so 
$O \subseteq \omega(y) \cup \alpha(y)$. 
Since $y \in \overline{O}  - O$, 
we obtain $y \in \omega(y) \cup \alpha(y)$. 
Then $y$ is not proper. 
The regularity of $y$ implies 
$y \notin \mathop{\mathrm{Sing}}(v)$ 
and so 
$y \in Q \setminus (\mathop{\mathrm{Sing}}(v) \sqcup \mathrm{P}) 
\subseteq \mathrm{LD} \sqcup \mathrm{E}$. 
Thus 
$\hat{O} \subseteq Q \setminus (\mathop{\mathrm{Sing}}(v) \sqcup \mathrm{P})$. 
On the other hand, 
we show that 
$\hat{O} \supseteq Q \setminus (\mathop{\mathrm{Sing}}(v) \sqcup \mathrm{P})$. 
Indeed, 
if there is exactly one quasi-minimal set $Q$, 
then 
$\overline{O_v(x)} = Q$ for any 
$x \in Q  \setminus (\mathop{\mathrm{Sing}}(v) \sqcup \mathrm{P}) 
\subseteq \mathrm{LD} \sqcup \mathrm{E}$. 
Thus we may assume that 
there are at least two quasi-minimal sets. 
The above Ma\v \i er work \cite{M} 
(cf. Remark 2 \cite{AZ}) impies that 
the genus of  $M$ is at least two.  
Note 
Cherry has proved that a quasi-minimal set contains 
a continuum of nontrivially recurrent orbits each of  which is dense
in the quasi-minimal set (Theorem VI\cite{C}).  
Therefore $Q$ contains 
nontrivially recurrent orbits. 
Fix a recurrent point $x \in (\mathrm{LD} \sqcup \mathrm{E}) \cap Q$ 
whose orbit closure is $Q$. 
For 
any point $y \in Q \setminus (\mathop{\mathrm{Sing}}(v) \sqcup \mathrm{P})$, 
we have 
$y \in \mathrm{LD} \sqcup \mathrm{E}$  
and so 
$y$ is weakly recurrent. 
By the above Cherry result, 
there is a recurrent point $z \in \overline{O_v(y)}$ 
whose orbit closure is $\overline{O_v(y)}$. 
By another Ma\v \i er theorem (cf. Theorem 4.2 \cite{AZ}) 
and its dual, 
we obtain 
$\omega(x) = \omega(z)$ 
and 
$\alpha(x) = \alpha(z)$.  
Thus 
$\overline{O} = 
Q = 
\overline{O_v(x)} = 
\overline{O_v(z)} = 
\overline{O_v(y)}$. 
This means  
$\hat{O} = Q \setminus (\mathop{\mathrm{Sing}}(v) \sqcup \mathrm{P})$.  
\end{proof}

%
We state a key lemma which is a relation between exceptional 
and proper orbits.
Recall that a subset $S$ of a surface $M$ is essential if 
some connected component of $S$ 
is neither null homotopic 
nor homotopic to a subset of the boundary $\partial M$.

\begin{lemma}\label{lem0bb}
$\overline{\mathop{\mathrm{Sing}}(v) \sqcup 
\mathop{\mathrm{Per}}(v) \sqcup \mathrm{LD}} 
\cap \mathrm{E}= \emptyset$ 
and 
$\overline{\mathop{\mathrm{Sing}}(v) \sqcup 
\mathop{\mathrm{Per}}(v) \sqcup \mathrm{E}} 
\cap \mathrm{LD}= \emptyset$.  
Moreover 
$\mathrm{E}\subset  \mathrm{int} \overline{\mathrm{P}}$. 
\end{lemma}

\begin{proof}
%
%
By taking a double covering of $M$ and the doubling of $M$ 
if necessary, 
we may assume that $v$ is transversally orientable 
and $M$ is closed and orientable. 
By the Ma\v \i er theorem, 
$\overline{\mathrm{E}}$
(resp. $\overline{\mathrm{LD}}$) 
consists of finitely many closures of 
exceptional (resp. locally dense) orbits.  
By Proposition \ref{lem001}, 
we have 
$\overline{\mathrm{LD}} \cap \mathrm{E} = \emptyset$ 
and 
$\mathrm{LD} \cap \overline{\mathrm{E}} = \emptyset$. 
Recall that 
the flow box theorem 
implies that 
the orbits of $v$ on a surface $M - \mathop{\mathrm{Sing}}(v)$ 
form a foliation $\mathcal{F}$. 
%
Moreover there is 
a transverse foliation $\mathcal{L}$ for this foliation $\mathcal{F}$ by 
Proposition 2.3.8 (p.18 \cite{HH}). 
%
%
%
%
We show 
$ \overline{\mathop{\mathrm{Per}}(v)} \cap (\mathrm{E} \sqcup \mathrm{LD}) =  \emptyset$.
%
Otherwise
there is a sequence of periodic orbits $O_i$ 
such that 
the closure $\overline{\cup_i O_i}$ 
contains 
a point $x \in \mathrm{E} \sqcup \mathrm{LD}$.  
By Lemma \ref{lem00}, 
we have  $\overline{O_v(x)} \cap \mathop{\mathrm{Per}}(v) = \emptyset$. 
We show that 
there is $K \in \mathbb{Z}_{>0}$ such that 
$O_k$ is contractible in $M_K$ 
for any $k > K$ 
where $M_K$ is the resulting closed surface 
of adding 
$2K$ disks to $M - (O_1 \sqcup \cdots \sqcup O_K)$. 
Indeed,  
we may assume that 
$O_1$ is 
essential 
by renumbering. 
Let $M_1$ be 
the resulting closed surface 
of adding two center disks to $M - O_1$. 
Then $g(M_1)< g(M)$, 
where $g(N)$ is the genus of a surface $N$. 
Since $M$ and $M_1$ are closed surfaces, 
by induction 
for essential closed curves 
at most $g(M)$ times, 
the assertion is followed. 
Since $M$ is normal, 
there are open disjoint neighbourhoods $U_x$ and $V$ of 
$\overline{O_v(x)}$ and $\cup_{i \leq K} O_i$ respectively.
Then 
there is 
a transverse arc $\gamma \subset U_x \cap l_x$ through $x$ 
which does not intersect $\cup_{i \leq K} O_i$, 
where $l_x \in \mathcal{L}$. 
Since $O_v(x)$ is exceptional or locally dense, 
we have 
$x$ is a  weakly recurrent point and so 
there is an arc $\gamma'$ in $O_v(x)$ 
whose boundaries are contained in $\gamma$. 
Since $M$ is normal, 
there is a neighbourhood of $\gamma \cup \gamma'$ 
which does not intersect $\cup_{i \leq K} O_i$. 
By the waterfall constraction (cf. Lemma 1.2, p.46\cite{ABZ}), 
we can construct 
a closed transversal $T$ for $v$ through $x$ 
which does not intersect $\cup_{i \leq K} O_i$. 
Let $v_K$ be a resulting flow on $M_K$ by adding center disks. 
Then $T$ is also a closed transversal for $v_K$. 
However there is $k > K$ such that 
$T$ intersects $O_k$ which is contractible in $M_K$. 
This is impossible. 
%
The closedness of $\mathop{\mathrm{Sing}}(v)$ 
implies the first assertion. 
Since $M = 
\mathop{\mathrm{Sing}}(v) \sqcup 
\mathop{\mathrm{Per}}(v) \sqcup 
\mathrm{LD}\sqcup \mathrm{P} \sqcup \mathrm{E}$, 
we obtain that 
$\mathrm{E}$ is contained in an open subset 
$M - \overline{\mathop{\mathrm{Sing}}(v) \sqcup 
\mathop{\mathrm{Per}}(v) \sqcup \mathrm{LD}} 
\subseteq \mathrm{P} \sqcup \mathrm{E}$.  
Therefore 
$\mathrm{E} \subseteq \mathrm{int}(\mathrm{P} \sqcup \mathrm{E})$. 
Since $\overline{\mathrm{E}}$ consists of finitely many closures of exceptional orbits, 
we have that 
${\mathrm{E}}$ is nowhere dense. 
This implies 
$\mathrm{E} \subset \overline{\mathrm{P}}$ 
and so  
$\mathrm{E} \subset \mathrm{int}(\mathrm{P} \sqcup \mathrm{E}) 
 \subset  \mathrm{int} \overline{\mathrm{P}}$. 
\end{proof}

From now on, 
we consider only non-wandering cases 
in this section.

\begin{lemma}\label{lem0aa}
Let $v$ a non-wandering flow on a compact surface $M$. 
Then 
$\mathop{\mathrm{Per}}(v)$ is open, 
$M = 
\mathop{\mathrm{Sing}}(v) \sqcup 
\mathop{\mathrm{Per}}(v) \sqcup 
\mathrm{LD}\sqcup \mathrm{P}$, 
and 
$\overline{\mathrm{LD}\sqcup \mathop{\mathrm{Per}}(v)}  
\supseteq M - \mathop{\mathrm{Sing}}(v)$. 
\end{lemma}

\begin{proof}
%
By taking a double covering of $M$ if necessary, 
we may assume that $v$ is transversally orientable. 
By Theorem III.2.12, III.2.15 \cite{BS}, 
the set of recurrence points is dense in $M$. 
By Lemma \ref{lem00}, 
there are no limit cycles. 
Write $U := M - \overline{\mathop{\mathrm{Sing}}(v) \sqcup 
\mathop{\mathrm{Per}}(v) \sqcup \mathrm{LD}}$.  
By Lemma \ref{lem0bb}, 
this $U \subseteq \mathrm{P} \sqcup \mathrm{E}$ is an open \nbd of 
$\mathrm{E}$. 
Since each point of $\mathrm{P}$ is not weakly recurrent, 
we have $\overline{\mathrm{E}} \supseteq U$. 
Since $\mathrm{E}$ is nowhere dense, 
we have $U$ is empty and so is $\mathrm{E}$. 
%
%
Hence 
$M = 
\mathop{\mathrm{Sing}}(v) \sqcup 
\mathop{\mathrm{Per}}(v) 
\sqcup 
\mathrm{LD}\sqcup \mathrm{P}$. 
Since there are no limit cycles, 
we obtain 
%
%
%
$\overline{\mathrm{LD} 
} \cap \mathop{\mathrm{Per}}(v) = \emptyset$.  
Fix a periodic orbit $O$. 
Then there is an annular  
\nbd $V$ of $O$ 
which is a finite union of flow boxes 
such that 
$V \cap 
(\mathop{\mathrm{Sing}}(v) \sqcup \overline{\mathrm{LD}} )
 = \emptyset$. 
Hence 
$V \subseteq 
\mathop{\mathrm{Per}}(v) \sqcup \mathrm{P}$.  
Since the set of recurrence points is dense in $M$, 
we have that 
$V \cap \mathop{\mathrm{Per}}(v)$ 
is dense in $V$. 
Therefore 
the holonomy of $O$ is 
identical 
and so 
$V \subseteq \mathop{\mathrm{Per}}(v)$. 
This means that 
$\mathop{\mathrm{Per}}(v)$ is open. 
%
%
%
%
For any $x \in \mathrm{P}$, 
since the regular weakly recurrent points form 
$\mathrm{Per}(v) \sqcup \mathrm{LD}$, 
by non-wandering property, 
each neighbourhood of $x$ meets 
$\mathrm{Per}(v) \sqcup \mathrm{LD}$ 
and so 
$\overline{\mathrm{LD}\sqcup \mathop{\mathrm{Per}}(v)} \supseteq \mathrm{P}$.  
\end{proof}

Now we state the characterization of non-wandering flows. 

\begin{theorem}\label{prop0c}
Let $v$ be a continuous flow on a compact surface $M$. 
Then 
$v$ is non-wandering 
if and only if 
$\overline{\mathrm{LD} \sqcup \mathop{\mathrm{Per}}(v)} \cup \mathop{\mathrm{Sing}}(v) 
= M$.  
In particular, 
if $v$ is non-wandering, 
then 
$\mathop{\mathrm{Per}}(v)$ 
is open 
and 
there are no exceptional orbits. 
\end{theorem}

\begin{proof}
%
Suppose that $v$ is non-wandering. 
By Lemma \ref{lem0aa}, 
we have 
$\overline{\mathrm{LD}\sqcup \mathop{\mathrm{Per}}(v)} \cup \mathop{\mathrm{Sing}}(v) 
= M$.  
Conversely, 
suppose that 
$\overline{\mathrm{LD}\sqcup \mathop{\mathrm{Per}}(v)} \cup \mathop{\mathrm{Sing}}(v) = M$.  
For any regular point $x$ of $M$, 
we have 
$x \in \overline{\mathrm{LD}\sqcup \mathop{\mathrm{Per}}(v)}$. 
This shows that $v$ is non-wandering. 
\end{proof}

We state the following charcterization of 
the union 
$\mathrm{P} \sqcup \mathop{\mathrm{Sing}}(v)$ of non-periodic proper orbits.

\begin{proposition}\label{lem32}
Let $v$ be a continuous non-wandering flow 
on a compact surface $M$. 
Then $\mathrm{P} \sqcup \mathop{\mathrm{Sing}}(v) = 
\{ x \in M \mid \omega(x) \cup \alpha(x) \subseteq \mathop{\mathrm{Sing}}(v) \}$. 
\end{proposition}

\begin{proof}
We may assume 
$M$ is connected. 
By taking a double covering of $M$ and the doubling of $M$ 
if necessary, 
we may assume that $v$ is transversally orientable 
and $M$ is closed and orientable. 
Obviously 
$\mathrm{P}  \sqcup \mathop{\mathrm{Sing}}(v) 
\supseteq \{ x \in M \mid \omega(x) \cup \alpha(x) \subseteq \mathop{\mathrm{Sing}}(v) \}$. 
Therefore 
it suffices to show this converse relation.
Fix a point $y \in \mathrm{P}$. 
By Lemma \ref{lem00}, 
the non-wandering property implies  
$\overline{O_v(y)} \cap \mathop{\mathrm{Per}}(v) = \emptyset$. 
We show 
$\overline{O_v(y)} \cap \mathrm{LD} = \emptyset$. 
Otherwise 
$\overline{O_v(y)}$ contains a locally dense orbit $O$. 
Since $\mathrm{int}\overline{O}$ is an open \nbd of $O$, 
we have  
$O_v(y) \cap \mathrm{int}\overline{O} \neq \emptyset$. 
The relation 
$y \in \overline{O} \subseteq \overline{O_v(y)}$ 
implies that 
$y$ is not proper, which is a contradiction. 
%
%
%
%
Since $\mathrm{E} = \emptyset$, 
we have 
$\overline{O_v(y)} \subset \mathrm{P} \sqcup \mathop{\mathrm{Sing}}(v)$. 
Suppose that 
$y \in \overline{\mathrm{LD}}$. 
Thus there is a recurrent point in $\mathrm{LD}$ 
whose orbit closure contains $O_v(y)$. 
Removing the singular points, 
Theorem 3.1\cite{Marz} implies 
$\omega(y) \cup \alpha(y) \subseteq \mathop{\mathrm{Sing}}(v)$. 
This means 
$y \in 
\{ x \in M \mid \omega(x) \cup \alpha(x) \subseteq \mathop{\mathrm{Sing}}(v) \}$. 
Otherwise 
$y \notin \overline{\mathrm{LD}}$. 
Suppose that 
each periodic orbit is null homotopic (i.e. non-essential). 
Fix a saturated \nbd $W$ of $y$. 
Then the intersection   
$W \cap \mathop{\mathrm{Per}}(v)$ 
is open dense in $W$.  
To apply Theorem 3.1\cite{Marz}, 
we need replace a small flow box near $O_v(y)$. 
Take a flow box $B \subset W$ 
which can be identified with $[-1, 1] \times [-1, 1]$ 
such that 
$O_v(y) \cap B = [-1, 1] \times \{ 0 \}$ 
and 
each orbit in $B$ is $[-1, 1] \times \{ a \}$ for some $a \in [-1, 1]$. 
Fix a small transverse arc $\gamma$ through $y$ in $B$. 
Let $\gamma_+$ be a connected component of $\gamma \setminus O$. 
Since each periodic orbit is null homotopic, 
each periodic orbit intersects $\gamma_+$ at most one point. 
This means that 
each point $ p:= (1, \varepsilon)$ in $\mathop{\mathrm{Per}}(v) \cap (\{ 1 \} \times [-1, 1])$ 
intersecting $\gamma_+$ 
goes back to a point $(-1, \varepsilon)$ in $B$ 
(i.e. $O_v(p) \cap B = \{ 1 \} \times \{ \varepsilon \}$). 
%
Therefore 
the flow $v|_B$ induces 
a homeomorphism 
from $\{-1 \} \times [-1, 1]$ 
to $\{ 1 \} \times [-1, 1]$ 
which can be identified with the identity mapping $1_{[-1, 1]}$ on $[-1, 1]$. 
Replacing $1_{[-1, 1]}$ 
with a homeomorphism $f$ which is contracting near $0$ 
(e.g. $f(x) = x^3$), 
we obtain the resulting continuous flow $w$ 
such that 
$O_v(y) = O_w(y)$ and 
$\mathop{\mathrm{Sing}}(w) = 
\mathop{\mathrm{Sing}}(v)$, 
modifying $v$ in $B$. 
We identify $O_v(y) \cap B = [-1, 1] \times \{ 0 \}$ 
with $\mathrm{dom} (f) = [-1, 1]$. 
By the Baire category theorem, 
the countable intersection 
$\bigcap_{n \in \Z} f^{n}([-1, 1] \cap \mathop{\mathrm{Per}}(v))$ 
is dense in $[-1, 1]$. 
Thus 
there is a point $z \in \mathop{\mathrm{Per}}(v)$ 
such that 
$O_w(z)$ is proper 
and 
$y \in \overline{O_w(z)} - O_w(z)$.  
Removing the singular points, 
Theorem 3.1\cite{Marz} for $w$ implies 
$\overline{O_w(y)} - O_w(y) = 
\alpha_w(y) \cup \omega_w(y) \subseteq \mathop{\mathrm{Sing}}(w)$. 
Since 
$O_v(y) = O_w(y)$ and 
$\mathop{\mathrm{Sing}}(w) = 
\mathop{\mathrm{Sing}}(v)$, 
we have 
$\omega(y) \cup \alpha(y) \subseteq \mathop{\mathrm{Sing}}(v)$. 
Suppose that 
there are essential periodic orbits. 
Let $O$ be an essential periodic orbit. 
Cutting this periodic orbit $O$, 
the remainder $M - O$ has 
new two boundaries. 
Adding two center disks to the two boundaries of $M - O$,  
we obtain 
a new closed surface $M'$ whose genus is less than one of $M$ 
and the resulting non-wandering flow $v'$ on $M'$ 
with 
$\mathop{\mathrm{Sing}}(v') = \mathop{\mathrm{Sing}}(v)$
such that 
$\alpha(y) = \alpha_{v'}(y)$ 
and 
$\omega(y) = \omega_{v'}(y)$, 
where 
$\alpha_{v'}(y)$ (resp. $\omega_{v'}(y)$) is 
the alpha (resp. omega) limit set of $y$ with respect to $v'$. 
By finite iterations, 
we obtain a new closed surface $M^{\dagger}$ 
and the resulting non-wandering flow $v^{\dagger}$ on $M^{\dagger}$ 
with 
$\mathop{\mathrm{Sing}}(v^{\dagger}) = \mathop{\mathrm{Sing}}(v)$
such that 
$\alpha(y) = \alpha_{v^{\dagger}}(y)$ 
and 
$\omega(y) = \omega_{v^{\dagger}}(y)$ 
such that 
each periodic orbit of $v^{\dagger}$ is null homotopic. 
This implies 
$\omega(y) \cup \alpha(y) \subseteq \mathop{\mathrm{Sing}}(v)$. 
\end{proof}

Taking a suspension of 
a non-wandering circle homeomorphism, 
we have the following statement.

\begin{corollary}\label{cor26}
Each non-wandering continuous homeomorphism 
on $\S^1$  
is topologically conjugate to 
a rotation. 
\end{corollary}

\begin{proof}
Let $v$ be the suspension of 
a non-wandering circle homeomorphism $f$. 
Then $v$ is a non-wandering continuous flow on $\T^2$.
Since $v$ is a suspension, 
 there are no singular points. 
 By Proposition \ref{lem32}, 
 we have 
 $\mathrm{P} = \emptyset$.  
Theorem \ref{prop0c} implies 
 $\mathop{\mathrm{Per}}(v) \sqcup \mathrm{LD} = \T^2$. 
 Since $\mathop{\mathrm{Per}}(v)$ 
is open, 
the complement 
$\mathrm{LD} = \T^2 - \mathop{\mathrm{Per}}(v)$ 
is closed. 
By Lemma \ref{lem0bb}, 
we obtain 
$\overline{\mathop{\mathrm{Per}}(v)} 
\cap \overline{\mathrm{LD}} 
= \overline{\mathop{\mathrm{Per}}(v)} 
\cap {\mathrm{LD}}  
= \emptyset$. 
This implies that 
the both of 
$\mathop{\mathrm{Per}}(v)$ 
and $\mathrm{LD}$ 
are open and closed. 
Since $\T^2$ is connected, 
we have that 
$v$ is pointwise periodic or 
minimal. 
Then $f$ is topologically conjugate to periodic or 
minimal. 
This means that 
$f$ is topologically conjugate to 
a rotation. 
\end{proof}

By the smoothing result \cite{G}, 
we have the following result.

\begin{corollary}\label{cor21}
Each non-wandering continuous flow on a compact surface $M$ 
is topologically equivalent to a $C^{\infty}$ flow. 
\end{corollary}

\begin{proof}
We may assume that 
$M$ is connected. 
By Theorem \ref{prop0c}, 
there are no exceptional orbits. 
By Proposition \ref{lem32}, 
the closure of each orbit in $\mathrm{P}$ contains 
singular points.  
We show that 
if the closure of locally dense orbit is minimal 
then it is the whole surface which is $\mathbb{T}^2$. 
Indeed, 
let $O$ be a locally dense orbit whose closure is minimal. 
Then $\overline{O}$ contains 
neither singular points
nor periodic points. 
Since the closure of an orbit in $\mathrm{P}$ contains 
singular points, 
the minimal set $\overline{O}$ consists of 
locally dense orbits. 
Since each point of $\overline{O}$ is contained 
in the interior of $\overline{O}$, 
this minimal set $\overline{O}$ is closed and open, and so  
is the whole surface  $M$ which is homeomorphic to $\mathbb{T}^2$.
Therefore each minimal set is 
either a closed orbit or $\mathbb{T}^2$. 
By the smoothing theorem \cite{G}, 
this continuous flow is topologically equivalent to a $C^{\infty}$ flow. 
\end{proof}

We state uniformity of $\mathop{\mathrm{Per}}(v)$ 
of a non-wandering continuous flow $v$. 

\begin{corollary}\label{cor29}
Let $v$ be a non-wandering continuous flow 
on a compact surface $M$.  
Then 
each connected component of $\mathop{\mathrm{Per}}(v)$ 
is either 
a connected component of $M$,   
an open annulus, or 
an open M\"obius band. 
If 
$\mathop{\mathrm{Sing}}(v)$ 
consists of finitely many contractible connected components, 
then $\mathrm{P}$ consists of finitely many separatrices,  
$\mathop{\mathrm{Sing}}(v) \sqcup \mathrm{P}$ is closed, 
and 
$\mathrm{LD}$ is open. 
\end{corollary}

\begin{proof} 
Since $\mathop{\mathrm{Per}}(v)$ is open, 
the flow box theorem 
implies that 
each periodic orbit has a saturated \nbd 
which is either an annulus or a M\"obius band. 
Notice that  
a union of one saturated M\"obius band and one saturated annulus with one intersection 
is a M\"obius band 
and that 
a union of two saturated annuli (resp. M\"obius bands) with one intersection 
is an annulus (resp. a Klein bottle).
Fix a connected component $U$ of $\mathop{\mathrm{Per}}(v)$.  
If $\partial U = \emptyset$, 
then $U$ is a connected component of $M$. 
If $\partial U \neq \emptyset$, 
then $U$ is 
either an annulus 
or a M\"obius band. 
Suppose that $\mathop{\mathrm{Sing}}(v)$ 
consists of finitely many contractible connected components. 
Collapsing each connected component of 
$\mathop{\mathrm{Sing}}(v)$ into a point, 
by Theorem 1\cite{RS}, 
the resulting space is homeomorphic to 
the original space $M$ 
and the resulting flow is non-wandering. 
Thus we may assume that 
$\mathop{\mathrm{Sing}}(v)$ is finite.  
By Theorem 3\cite{CGL}, 
each singular point is 
either a center or a multi-saddle. 
By Proposition \ref{lem32}, 
each orbit in $\mathrm{P}$ is 
a separatrix. 
Therefore 
$\mathrm{P}$ consists of finitely many separatrices and so   
$\mathop{\mathrm{Sing}}(v) \sqcup \mathrm{P}$ is closed.   
Since $\overline{\mathop{\mathrm{Per}}(v)} \cap \mathrm{LD} = \emptyset$, 
the complement 
$\mathrm{LD} = M - (\mathop{\mathrm{Sing}}(v) \sqcup \mathrm{P} \sqcup \mathop{\mathrm{Per}}(v))$ is open. 
\end{proof}

We show that 
$\mathrm{LD}$ is not open in general.  
In fact, 
we construct 
a non-wandering flow on $\T^2$ 
such that 
$\mathrm{P}$ and 
$\mathrm{LD}$ are dense as follows. 

\begin{example}\label{ex1}
Consider an irrational rotation $v$ on $\T^2$. 
Fix any points $p \in \T^2$. 
Let  $( t_i )_{i \in \Z}$ of $\R$ be 
a sequence such that 
$\lim_{i \to \infty} t_i = \infty$, 
$\lim_{i \to - \infty} t_i = - \infty$, 
and 
$\lim_{i \to \infty} p_i 
= \lim_{i \to - \infty} p_i 
= p_0$, 
where $p_i := v(t_i, p_0) \in O_v(p)$.  
Using dump functions, 
replace $O_v(p)$ 
with a union of countably many singular points $p_i$ ($i \in \Z$) 
and countably many proper orbits. 
%
Let $v'$ be the resulting vector field. 
For any point $x \in \T^2 - O(p)$, 
we have $O_v(x) = O_{v'}(x)$. 
Moreover $O_v(p) - \mathop{\mathrm{Sing}}(v') = \mathrm{P}(v')$.  
By construction, 
$\mathrm{LD}$ is not open but 
$\overline{\mathrm{P}} = \overline{\mathrm{LD}} = \T^2$. 
This example also shows that 
the finiteness condition in 
Corollary \ref{cor29} 
is necessary. 
\end{example}

In \cite{H2}, 
a foliation $\F$  on a manifold $M$ is 
said to be ``rare species'' 
if 
either 
all leaves are exceptional 
or 
$\F$ has at least two of three type (i.e. proper, locally dense, exceptional) 
and the union of leaves of each type is dense. 
The author 
have constructed 
some kind of   
codimension one ``rare species''  foliations on compact $3$-manifolds.
Analogically we call that 
a flow $v$ on $M$ is 
``rare species''
if 
either 
$\mathrm{E} = M$ 
or 
$v$ has at least two of three type (i.e. proper, locally dense, exceptional) 
and the union of orbits of each type is dense. 
Note 
that ``rare species'' surface flows 
are non-wandering and so 
have no exceptional orbits. 
In contrast to 
foliations, 
this implies 
there is only one possible kind of 
``rare species'' flows on compact surfaces.  
Notice that 
the above example is 
a ``rare species'' flows on $\mathbb{T}^2$. 
Therefore, 
we obtain the following statement.

\begin{proposition}\label{prop}
There are smooth ``rare species'' flows on $\mathbb{T}^2$.   
On the other hand, 
for each ``rare species'' flow $v$ on a compact surface $M$, we have 
$M = \mathop{\mathrm{Sing}}(v) \sqcup \mathrm{LD} \sqcup \mathrm{P}$ 
and 
$\overline{\mathrm{LD}} = 
\overline{\mathrm{P}} = M$.   
\end{proposition}

Note that 
all known codimension one foliations on 
compact manifolds are not $C^1$ but continuous. 
In contrast to codimension one foliations on 
compact manifolds,  
each ``rare species'' flow $v$ on a compact surface 
is topologically equivalent 
to a $C^{\infty}$ flow. 
Moreover, 
we obtain following examples. 

\begin{corollary}\label{cor}
There are codimension one ``rare species'' foliations 
on open surfaces contained in $\R^2$.  
In particular, 
the union of locally dense leaves 
is not open. 
\end{corollary}

\begin{proof}
Fix a minimal codimension one foliation $\F$ 
on open surfaces contained in $\R^2$  
(e.g. a foliation in Theorem \cite{F}). 
Replace a locally dense leaf $L$ 
into countably many singular points and 
proper leaves connecting two singular points,  
and  
remove the singular points as above construction for flows. 
Then we can obtain a desired foliation. 
\end{proof}

\section{Applications of this characterization}

Recall that 
$v$ is topologically transitive if 
it has a dense orbit.  
Closed orbits are singular points or periodic orbits.
A subset is said to be co-connected if 
the complement of it is connected. 
We have a following characterization of 
transitivity for surface flows. 

\begin{theorem}\label{th0}
Let $v$ be a continuous flow on a compact surface $M$. 
Then 
the following are equivalent: 
\\
1. 
$v$ is topologically transitive. 
\\
2.  
$v$ is non-wandering 
such that 
$\mathrm{P} \sqcup \mathop{\mathrm{Sing}}(v)$ is co-connected 
and 
$ \mathop{\mathrm{int}}(\mathop{\mathrm{Per}}(v) 
\sqcup \mathop{\mathrm{Sing}}(v)) = \emptyset$.  
\\
3. 
$v$ is non-wandering 
such that 
the set of regular weakly recurrent points is connected 
and 
the interior of 
the union of closed orbits is empty.  
\\
4. 
$\mathrm{P} \sqcup \mathop{\mathrm{Sing}}(v)$ is co-connected 
and 
$\mathop{\mathrm{int}} \mathop{\mathrm{Sing}}(v) = 
\mathop{\mathrm{int}} \mathop{\mathrm{Per}}(v) = 
\mathop{\mathrm{int}} \mathrm{P} = 
\emptyset$.  

In each case, 
$M 
= \mathop{\mathrm{Sing}}(v) \sqcup 
\mathrm{P} \sqcup \mathrm{LD} 
= \overline{\mathrm{LD}}$ 
and 
each locally dense orbit is dense. 
\end{theorem}

\begin{proof}
Obviously, 
the conditions 2 and 3 are equivalent. 
Since $\mathrm{P}$ is the complement of 
the set of the weakly recurrent points, 
$v$ is non-wandering if 
and only if $\mathop{\mathrm{int}} \mathrm{P} = \emptyset$.  
Lemma \ref{lem0aa} implies that 
the conditions 2 and 4 are equivalent. 
Suppose that 
$v$ is topologically transitive. 
Then $v$ is non-wandering. 
The openness of  $\mathop{\mathrm{Per}}(v)$ 
implies $\mathop{\mathrm{Per}}(v) = \emptyset$. 
%
By transitivity, 
there is a dense orbit $O$ 
and so $\overline{\mathrm{LD}} = M$. 
Since $O$ is connected 
and $O \subseteq \mathrm{LD} 
\subseteq \overline{O} = M$, 
we have that 
$M - (\mathrm{P} \sqcup \mathop{\mathrm{Sing}}(v)) = \mathrm{LD}$ is connected. 
Conversely, 
suppose 
$v$ is non-wandering 
such that 
$ \mathop{\mathrm{int}}(\mathop{\mathrm{Per}}(v) 
\sqcup \mathop{\mathrm{Sing}}(v)) = \emptyset$ 
and 
$M - (\mathrm{P} 
\sqcup \mathop{\mathrm{Sing}}(v))$ is connected. 
Put $U = M - (\mathrm{P} 
\sqcup \mathop{\mathrm{Sing}}(v))$. 
Since $\mathop{\mathrm{Per}}(v)$ is open, 
we have 
$\mathop{\mathrm{Per}}(v)  = \emptyset$ 
and so 
$\mathop{\mathrm{int}}\mathop{\mathrm{Sing}}(v) = \emptyset$.  
Since $v$ is non-wandering, 
we have 
$\overline{\mathrm{LD}} = M = 
\mathop{\mathrm{Sing}}(v) \sqcup 
\mathrm{LD} \sqcup \mathrm{P}$. 
Thus $\mathrm{LD} = U$ is connected. 
Fix an locally dense orbit $O$. 
For any $z \in \hat{O}$, 
Proposition \ref{lem32} implies 
$z \in \mathrm{LD}$ 
and 
so 
$z \in \mathrm{int}\overline{O_v(z)} = \mathrm{int}\overline{O}$. 
Hence 
$\hat{O} \subseteq \mathrm{int} \overline{O}$ 
and so 
$\hat{O} = 
\overline{O} \cap \mathrm{LD} = 
\mathrm{int} \overline{O} \cap \mathrm{LD}$. 
This means that 
$\hat{O}$ is closed and open in $\mathrm{LD}$ 
and so 
$\overline{O} \cap \mathrm{LD} = \mathrm{LD}$. 
Therefore  
$M = 
\overline{\mathrm{LD}} = 
\overline{O}$. 
%
\end{proof}

Since an area-preserving flow is non-wandering, 
Proposition \ref{lem32} implies 
the following statement which is a generalization of Theorem A \cite{MS}.

\begin{corollary}\label{cor}
Let $v$ be a continuous flow on a compact surface $M$. 
Then 
$v$ is topologically transitive 
if and only if 
$v$ is non-wandering 
such that 
$ \mathop{\mathrm{int}}(\mathop{\mathrm{Per}}(v) 
\sqcup \mathop{\mathrm{Sing}}(v)) = \emptyset$ 
and 
$\{ x \in M \mid \omega(x) \cup \alpha(x) \subseteq \mathop{\mathrm{Sing}}(v) \}$ is co-connected. 
\end{corollary}

\end{document}